# STATISTICAL INFERENCE FOR TIME-INHOMOGENEOUS VOLATILITY MODELS

By Danilo Mercurio[1] and Vladimir Spokoiny

*Humboldt-Universität zu Berlin and Weierstrass-Institute*

This paper offers a new approach for estimating and forecasting the volatility of financial time series. No assumption is made about the parametric form of the processes. On the contrary, we only suppose that the volatility can be approximated by a constant over some interval. In such a framework, the main problem consists of filtering this *interval of time homogeneity*; then the estimate of the volatility can be simply obtained by local averaging. We construct a *locally adaptive volatility estimate* (LAVE) which can perform this task and investigate it both from the theoretical point of view and through Monte Carlo simulations. Finally, the LAVE procedure is applied to a data set of nine exchange rates and a comparison with a standard GARCH model is also provided. Both models appear to be capable of explaining many of the features of the data; nevertheless, the new approach seems to be superior to the GARCH method as far as the out-of-sample results are concerned.

**1. Introduction.** The aim of this paper is to offer a new perspective for the estimation and forecasting of the volatility of financial asset returns such as stocks and exchange rate returns.

A remarkable amount of statistical research is devoted to financial time series, in particular, to the volatility of asset returns, where the term volatility indicates a measure of dispersion, usually the variance or the standard deviation. The interest in this topic is motivated by the needs of the financial industry, which regards volatility as one of the main reference numbers for risk management and derivative pricing.

Actually, asset returns time series display very peculiar stylized facts, which are connected with their second moments. Graphically, they look like

Received June 2000; revised March 2003.

[1]Supported by the Deutsche Forschungsgemeinschaft, Graduiertenkolleg für Angewandte Mikroökonomik at Humboldt University.

*AMS 2000 subject classifications.* Primary 62M10; secondary 62P20.

*Key words and phrases.* Stochastic volatility model, adaptive estimation, local homogeneity.

white noise, where periods of high and low volatility seem to alternate. Their density has fat tails if compared to that of a normal random variable, and they show significantly positive and highly persistent autocorrelation of the absolute returns, meaning that large (resp. small) absolute returns are likely to be followed by large (resp. small) absolute returns. Typical examples can be seen in Section 6, and further details on this topic can be found in Taylor (1986). Therefore, a white-noise process with time-varying variance is usually taken to model such features. Let $S_t$ denote the observed asset process. Then the corresponding (log) returns $R_t = \log(S_t/S_{t-1})$ follow the heteroscedastic model

$$R_t = \sigma_t \xi_t,$$

where $\xi_t$ are standard Gaussian independent innovations and $\sigma_t$ is a time-varying *volatility* coefficient. It is often assumed that $\sigma_t$ is measurable w.r.t. the $\sigma$-field generated by the preceding returns $R_1, \ldots, R_{t-1}$. For modeling this volatility process, parametric assumptions are usually used. The main model classes are the ARCH and GARCH family [Engle (1995)] and the stochastic volatility models [Harvey, Ruiz and Shephard (1994)]. A large number of papers has followed the first publications on this topic, and the original models have been extended in order to provide better explanations. For example, models which take into account asymmetries in volatility have been proposed, such as EGARCH [Nelson (1991)], QGARCH [Sentana (1995)] and GJR [Glosten, Jagannathan and Runkle (1993)]; furthermore, the research on integrated processes has produced integrated [Engle and Bollerslev (1986)] and fractal integrated versions of the GARCH model.

The availability of very large samples of financial data has made it possible to construct models which display quite complicated parameterizations in order to explain all the observed stylized facts. Obviously, these models rely on the assumption that the parametric structure of the process remains constant through the whole sample. This is a nontrivial and possibly dangerous assumption, in particular, as far as forecasting is concerned [Clements and Hendry (1998)]. Furthermore, checking for parameter instability becomes quite difficult if the model is nonlinear and/or the number of parameters is large. Thus, those characteristics of the returns, which are often explained by the long memory and (fractal) integrated nature of the volatility process, could also depend on the parameters being time varying.

In this paper we propose another approach focusing on a very simple model but with a possibility for model parameters to depend on time. This means that the model is regularly checked and adapted to the data. No assumption is made about the parametric structure of the volatility process. We only suppose that it can be locally approximated by a constant; that is, for every time point $\tau$ there exists a past interval $[\tau - m, \tau]$ where the volatility $\sigma_t$ did not vary much. This interval is referred to as the *interval of time homogeneity*. An algorithm is proposed for data-driven estimation of the interval of time homogeneity, after which the estimate of the volatility can be simply obtained by averaging.

Our approach is similar to varying-coefficient modeling from Fan and Zhang (1999); see also Cai, Fan and Li (2000) and Cai, Fan and Yao (2000). Fan, Jiang, Zhang and Zhou (2003) discussed applications of this method to stock price volatility modeling. The proposed procedure is based on the assumption that the model parameters smoothly vary with time and can be locally approximated by a linear function of time. This approach has the drawback of not allowing one to incorporate structural breaks into the model.

Change point modeling with applications to financial time series was considered in Mikosch and Starica (2000). Kitagawa (1987) applied non-Gaussian random walk modeling with heavy tails as the prior for the piecewise constant mean for one-step-ahead prediction of nonstationary time series. However, the aforementioned approaches require some essential amount of prior information about the frequency of change points and their size.

The LAVE approach proposed in this article does not assume smooth or piecewise constant structure of the underlying process and does not require any prior information. The procedure proposed below in Section 3 focuses on adaptive choice of the interval of homogeneity that allows one to proceed in a unified way with smoothly varying coefficient models and change point models.

The proposed approach attempts to describe the *local* dynamic of the volatility process, and it is particularly appealing for short-term forecasting purposes which is an important building block, for example, in value-at-risk and portfolio hedging problems or backtesting [Härdle and Stahl (1999)].

The remainder of the paper is organized as follows. Section 2 introduces the adaptive modeling procedure. Then some theoretical properties are discussed in the general situation and for a change point model. A simulation study illustrates the performance of the new methodology with respect to the change point model. The question of selecting the smoothing parameters is also addressed and some solutions are proposed. Finally, the procedure is applied to a set of nine exchange rates and it appears to be highly competitive with standard GARCH$(1,1)$, which is used as a benchmark model. Mathematical proofs are given in Section 8.

**2. Modeling volatility via power transformation.** Let $S_t$ be an observed asset process in discrete time, $t = 1, 2, \ldots, \tau$ and $R_t$ are the corresponding returns: $R_t = \log(S_t/S_{t-1})$. We model this process via the *conditional heteroscedasticity* assumption

$$R_t = \sigma_t \xi_t, \tag{2.1}$$

where $\xi_t$, $t \geq 1$, is a sequence of independent standard Gaussian random variables and $\sigma_t$ is the *volatility* process which is in general a predictable random process, that is, $\sigma_t \sim \mathcal{F}_{t-1}$ with $\mathcal{F}_{t-1} = \sigma(R_1, \ldots, R_{t-1})$ (the $\sigma$-field generated by the first $t-1$ observations).

A *time-homogeneous* (*time-homoscedastic*) model means that $\sigma_t$ is a constant. The process $S_t$ is then a geometric Brownian motion observed at discrete time moments. The assumption of time homogeneity is too restrictive in practical applications, and it does not allow one to fit real data very well. In this paper, we consider an approach based on the *local time homogeneity*, which means that for every time moment $\tau$ there exists a time interval $[\tau - m, \tau]$ where the volatility process $\sigma_t$ is nearly constant. Under such a modeling, the main intention is both to describe the interval of homogeneity and to estimate the corresponding value $\sigma_\tau$ which can then be used for one-step forecasting and the like.

2.1. *Data transformation.* The model equation (2.1) links the target volatility process $\sigma_t$ with the observations $R_t$ via the multiplicative errors $\xi_t$. The classical well-developed regression approach relies on the assumption of additive errors which can then be smoothed out by some kind of averaging. A natural and widespread method of transforming equation (2.1) into a regression-like equation is to apply the log function to both its sides squared:

$$\log R_t^2 = \log \sigma_t^2 + \log \xi_t^2, \tag{2.2}$$

which can be rewritten in the form

$$\log R_t^2 = \log \sigma_t^2 + C + v\zeta_t,$$

with $C = \mathbf{E} \log \xi_t^2$, $v^2 = \operatorname{Var} \log \xi_t^2$ and $\zeta_t = v^{-1}(\log \xi_t^2 - C)$; see, for example, Gouriéroux (1997). This is a usual regression equation with the "response" $Y_t = \log R_t^2$, target regression function $f(t) = \log \sigma_t^2 + C$ and homogeneous "noise" $v\zeta_t$.

The main problem with this approach is due to the distribution of the errors $\zeta_t$, which is highly skewed and gives very high weights to the small values of the errors $\xi_t$. In particular, this leads to a serious problem with missing data which are typically modeled equal to previous values providing $R_t = 0$.

Another possibility is based on power transformation [see Carroll and Ruppert (1988)] which also leads to a regression with additive noise and this noise is much closer to a Gaussian one. Due to (2.1), the random variable $R_t$ is conditionally on $\mathcal{F}_{t-1}$ Gaussian and

$$\mathbf{E}(R_t^2 | \mathcal{F}_{t-1}) = \sigma_t^2.$$

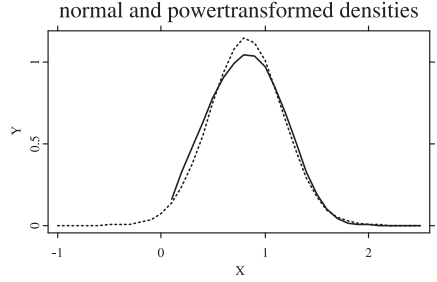

Fig. 1. *Density of $p_{1/2}(x)$ (straight line) and the standard normal density (dotted line).*

Similarly, for every $\gamma > 0$,

$$\mathbf{E}(|R_t|^\gamma | \mathcal{F}_{t-1}) = \sigma_t^\gamma \mathbf{E}(|\xi|^\gamma | \mathcal{F}_{t-1}) = C_\gamma \sigma_t^\gamma,$$

$$\mathbf{E}(|R_t|^\gamma - C_\gamma \sigma_t^\gamma | \mathcal{F}_{t-1})^2 = \sigma_t^{2\gamma} \mathbf{E}(|\xi|^\gamma - C_\gamma)^2 = \sigma_t^{2\gamma} D_\gamma^2,$$

where $\xi$ denotes a standard Gaussian r.v., $C_\gamma = \mathbf{E}|\xi|^\gamma$ and $D_\gamma^2 = \operatorname{Var}|\xi|^\gamma$. Therefore, the process $|R_t|^\gamma$ allows for the representation

$$(2.3) \qquad |R_t|^\gamma = C_\gamma \sigma_t^\gamma + D_\gamma \sigma_t^\gamma \zeta_t,$$

where $\zeta_t$ is equal to $(|\xi|^\gamma - C_\gamma)/D_\gamma$. Note that the problem of estimating $\sigma_t$ is in some sense equivalent to the problem of estimating $\theta_t = C_\gamma \sigma_t^\gamma$, which is the conditional mean of the transformed process $|R_t|^\gamma$. This is already a kind of heteroscedastic regression problem with additive errors $D_\gamma \sigma_t^\gamma \zeta_t$ satisfying

$$\mathbf{E}(D_\gamma \sigma_t^\gamma \zeta_t | \mathcal{F}_{t-1}) = 0,$$

$$\mathbf{E}(D_\gamma^2 \sigma_t^{2\gamma} \zeta_t^2 | \mathcal{F}_{t-1}) = D_\gamma^2 \sigma_t^{2\gamma}.$$

A natural choice of the parameter $\gamma$ is $\gamma = 2$ providing nearly efficient variance estimation under homogeneity. For $\gamma = 2$ one has $C_\gamma = 1$ and $D_\gamma^2 = 2$. Note, however, that the distribution of the "errors" $\zeta_t = (|\xi_t|^\gamma - C_\gamma)/D_\gamma$ is still heavy tailed and highly skewed, which results in a low sensitivity of the method in an inhomogeneous situation. The other important cases are $\gamma = 1$ and $\gamma = 1/2$. A minimization of the skewness $\mathbf{E}\zeta_\gamma^3$ and the fat $\mathbf{E}\zeta_\gamma^4 - 3$ with respect to $\gamma$ leads to the choice $\gamma \approx 1/2$. The corresponding density $p_{1/2}(x)$ of $\zeta_{1/2}$ together with the standard normal density $\phi(x)$ is plotted in Figure 1. Our numerical results are also in favor of the choice $\gamma = 1/2$; see Section 5.

**3. Adaptive estimation under local time homogeneity.** Here we describe one approach to volatility modeling based on the assumption of local time homogeneity starting from the preliminary heuristic discussion. The assumption of local time homogeneity means that the function $\sigma_t$ is nearly constant within an interval $I = [\tau - m, \tau]$, and the process $R_t$ follows the regression-like equation (2.3) with the constant trend $\theta_I = C_\gamma \sigma_I^\gamma$ which can be estimated by averaging over this interval $I$:

$$\tilde{\theta}_I = \frac{1}{|I|} \sum_{t \in I} |R_t|^\gamma. \tag{3.1}$$

For the particular case $\gamma = 2$ the estimate $\tilde{\theta}_I$ coincides with the local maximum likelihood estimator (MLE) of the volatility $\sigma_t^2$ considered in Fan, Jiang, Zhang and Zhou (2003). As discussed in the previous section, a smaller value of $\gamma$ might be preferred for improving the stability of the method. Similarly to Fan, Jiang, Zhang and Zhou (2003), one can also incorporate the one-sided kernel weighting to this estimator.

By (2.3)

$$\tilde{\theta}_I = \frac{C_\gamma}{|I|} \sum_{t \in I} \sigma_t^\gamma + \frac{D_\gamma}{|I|} \sum_{t \in I} \sigma_t^\gamma \zeta_t = \frac{1}{|I|} \sum_{t \in I} \theta_t + \frac{s_\gamma}{|I|} \sum_{t \in I} \theta_t \zeta_t, \tag{3.2}$$

with $s_\gamma = D_\gamma / C_\gamma$ so that

$$\mathbf{E}\tilde{\theta}_I = \mathbf{E}\frac{1}{|I|} \sum_{t \in I} \theta_t, \tag{3.3}$$

$$\frac{s_\gamma^2}{|I|^2} \mathbf{E}\left(\sum_{t \in I} \theta_t \zeta_t\right)^2 = \frac{s_\gamma^2}{|I|^2} \mathbf{E}\sum_{t \in I} \theta_t^2. \tag{3.4}$$

3.1. *Some properties of the estimate $\tilde{\theta}_I$.* Due to our assumption of local homogeneity, the process $\theta_t$ is close to $\theta_\tau$ for all $t \in I$. Define also

$$\Delta_I = \sup_{t \in I} |\theta_t - \theta_\tau| \quad \text{and} \quad v_I^2 = \frac{s_\gamma^2}{|I|^2} \sum_{t \in I} \theta_t^2.$$

The value of $\Delta_I$ measures the departure from homogeneity within the interval $I$, and it can be regarded as an upper bound of the "bias" of the estimate $\tilde{\theta}_I$. The value of $v_I^2$, because of (3.4), will be referred as the "conditional variance" of the estimate $\tilde{\theta}_I$. The next theorem provides a probability bound for the estimation error, that is, the deviation of $\tilde{\theta}_I$ from the present value of the volatility $\theta_\tau$ in terms of $\Delta_I$ and $v_I$.

THEOREM 3.1. *Let the volatility coefficient $\sigma_t$ satisfy the condition*

$$b \leq \sigma_t^\gamma \leq bB, \tag{3.5}$$

*with some positive constants $b, B$. Then there exists $a_\gamma > 0$ such that, for every $\lambda \geq 1$,*

$$\mathbf{P}(|\tilde{\theta}_I - \theta_\tau| > \Delta_I + \lambda v_I) \leq 4\sqrt{e}a_\gamma^{-1}\lambda(1 + \log B)e^{-\lambda^2/(2a_\gamma)}.$$

REMARK 3.1. This result can be slightly refined for the special case when the volatility process $\sigma_t$ for $t \in I$ is deterministic or (conditionally) independent of the observations $R_t$ preceding $I$. Namely, in such a situation the factor $4\sqrt{e}a_\gamma^{-1}\lambda(1 + \log B)$ in the bound can be replaced by 2:

$$\mathbf{P}(|\tilde{\theta}_I - \theta_\tau| > \Delta_I + \lambda v_I) \leq 2e^{-\lambda^2/(2a_\gamma)}.$$

A similar remark applies to all the results that follow.

The result of this theorem bounds the loss of the estimate $\tilde{\theta}_I$ via the value $\Delta_I$ and the conditional standard deviation $v_I$. Under homogeneity $\Delta_I \equiv 0$ and the error of estimation is of order $v_I$. Unfortunately, $v_I$ depends, in turn, on the target process $\theta_t$. One would be interested in another bound which does not involve the unknown function $\theta_t$. Namely, using (3.4) and assuming $\Delta_I$ small, one may replace the conditional standard deviation $v_I$ by its estimate

$$\tilde{v}_I = s_\gamma \tilde{\theta}_I |I|^{-1/2}.$$

THEOREM 3.2. *Let $R_1, \ldots, R_\tau$ obey (2.1) and let (3.5) hold true. Then, for the estimate $\tilde{\theta}_I$ of $\theta_\tau$ for every $D \geq 0$ and $\lambda \geq 1$,*

$$\mathbf{P}(|\tilde{\theta}_I - \theta_\tau| > \lambda' \tilde{v}_I, \Delta_I/v_I \leq D) \leq 4\sqrt{e}\lambda(1 + \log B)e^{-\lambda^2/(2a_\gamma)},$$

*where $\lambda'$ solves*

$$\lambda + D = \lambda'/(1 + \lambda' s_\gamma |I|^{-1/2}).$$

3.2. *Adaptive choice of the interval of homogeneity.* Given observations $R_1, \ldots, R_\tau$ following the time-inhomogeneous model (2.1), we aim to estimate the current value of the parameter $\theta_\tau$ using the estimate $\tilde{\theta}_I$ with a properly selected time interval $I$ of the form $[\tau - m, \tau]$ to minimize the corresponding estimation error. Below we discuss one approach which goes back to the idea of pointwise adaptive estimation; see Lepski (1990), Lepski and Spokoiny (1997) and Spokoiny (1998). The idea of the method can be explained as follows. Suppose $I$ is an interval candidate; that is, we expect time homogeneity in $I$ and, hence, in every subinterval of $I$. This particularly implies that the value $\Delta_I$ is small and similarly for all $\Delta_J$, $J \subset I$, and that the mean values of the $\theta_t$ over $I$ and over $J$ nearly coincide. Our adaptive procedure roughly means the choice of the largest possible interval $I$ such that the hypothesis that the value $\theta_t$ is a constant within $I$ is not rejected. For testing this hypothesis, we consider the family of subintervals of $I$ of the form $J = [\tau - m', \tau]$ with $m' < m$ and for every such subinterval $J$ compare two different estimates: one is based on the observations from $J$, and the other one is calculated from the complement $I \setminus J = [\tau - m, \tau - m'[$. Theorems 3.1 and 3.2 can be used to bound the difference $\tilde{\theta}_J - \tilde{\theta}_{I\setminus J}$ under homogeneity within $I$. Indeed, the conditional variance of $\tilde{\theta}_{I\setminus J} - \tilde{\theta}_J$ is $v_{I\setminus J}^2 + v_J^2$ and can be estimated by $\tilde{v}_{I\setminus J}^2 + \tilde{v}_J^2$. Thus, with high probability it holds that

$$|\tilde{\theta}_{I\setminus J} - \tilde{\theta}_J| \leq \lambda \sqrt{\tilde{v}_{I\setminus J}^2 + \tilde{v}_J^2},$$

provided that $\lambda$ is sufficiently large. Therefore, if there exists a testing interval $J \subset I$ such that the quantity $|\tilde{\theta}_{I\setminus J} - \tilde{\theta}_J|$ is significantly positive, then we reject the hypothesis of homogeneity for the interval $I$. Finally, our adaptive estimate corresponds to the largest interval $I$ such that the hypothesis of homogeneity is not rejected for $I$ itself and all smaller considered intervals.

Now we present a formal description. Suppose a family $\mathcal{I}$ of interval candidates $I$ is fixed. Each of them is of the form $I = [\tau - m, \tau]$, $m \in \mathbb{N}$, so that the set $\mathcal{I}$ is ordered due to $m$. With every such interval, we associate the estimate $\tilde{\theta}_I$ of $\theta_\tau$ and the corresponding estimate $\tilde{v}_I$ of the conditional standard deviations $v_I$.

Next, for every interval $I$ from $\mathcal{I}$ we assume there is a set $\mathcal{J}(I)$ of testing subintervals $J$ [one example of these sets $\mathcal{I}$ and $\mathcal{J}(I)$ is given in Section 6]. For every $J \in \mathcal{J}(I)$ we construct the corresponding estimate $\tilde{\theta}_J$ (resp. $\tilde{\theta}_{I\setminus J}$) from the observations $Y_t = |R_t|^\gamma$ for $t \in J$ (resp. for $t \in I \setminus J$) according to (3.1) and compute $\tilde{v}_J$ (resp. $\tilde{v}_{I\setminus J}$).

Now, with a constant $\lambda$, define the adaptive choice of the interval of homogeneity by the following iterative procedure:

*Initialization.*   Select the smallest interval in $\mathcal{I}$.

*Iteration.*   Select the next interval $I$ in $\mathcal{I}$ and calculate the corresponding estimate $\tilde{\theta}_I$ and the estimated conditional standard deviation $\tilde{v}_I$.

*Testing homogeneity.*   Reject $I$ if there exists one $J \in \mathcal{J}(I)$ such that

$$(3.6) \qquad |\tilde{\theta}_{I\setminus J} - \tilde{\theta}_J| > \lambda \sqrt{\tilde{v}_{I\setminus J}^2 + \tilde{v}_J^2}.$$

*Loop.*   If $I$ is not rejected, then continue with the iteration step by choosing a larger interval. Otherwise, set $\hat{I} =$ "the latest nonrejected $I$."

The *locally adaptive volatility estimate* (LAVE) $\hat{\theta}_\tau$ of $\theta_\tau$ is defined by applying this selected interval $\hat{I}$:

$$\hat{\theta}_\tau = \tilde{\theta}_{\hat{I}}.$$

The next section discusses the theoretical properties of the LAVE algorithm in a general framework, while Section 6 gives a concrete example for the choice of the sets $\mathcal{I}$, $\mathcal{J}(I)$ and the parameter $\lambda$. This choice is then applied to simulated and real data.

**4. Theoretical properties.** In this section we collect some results describing the quality of the proposed adaptive procedure.

4.1. *Accuracy of the adaptive estimate.* Let $\hat{I}$ be the interval selected by our adaptive procedure. We aim to show that our adaptive choice is up to some constant factor in the losses as good as the "ideal" choice $\mathbb{I}$ that may utilize the knowledge of the volatility process $\sigma_t$. This "ideal" choice can be defined by balancing the accuracy of approximating the underlying process $\theta_t$ (which is controlled by $\Delta_I$) and the stochastic error controlled by the stochastic standard deviation $v_I$. By definition, $v_I = s_\gamma |I|^{-1}(\sum_{t\in I}\theta_t^2)^{1/2}$ so that $v_I$ typically decreases when $|I|$ increases. For simplicity of notation we shall suppose further that $v_I \leq v_J$ for $J \subset I$.

We do not give a formal definition of an "ideal" choice of the interval $I$ since there is no one universally optimal choice even if the process $\theta_t$ is known. Instead, we consider a family of all "good" intervals $\mathbb{I}$ such that the variability of the process $\theta_t$ inside $\mathbb{I}$ is not too large compared to the conditional stochastic deviation $v_{\mathbb{I}}$. This, due to Theorem 3.1, allows us to bound with high probability the losses of the "ideal" estimate $\tilde{\theta}_{\mathbb{I}}$ by $(D+\lambda)v_{\mathbb{I}}$ provided that $\Delta_{\mathbb{I}}/v_{\mathbb{I}} \leq D$ and $\lambda$ is sufficiently large. A similar property should hold for all smaller intervals $I \subset \mathbb{I}$. Hence, it is natural to quantify the quality of the interval $\mathbb{I}$ by

$$\delta_{\mathbb{I}} = \sup_{I \in \mathcal{I}\,:\,I \subseteq \mathbb{I}} \Delta_I/v_I.$$

The next assertion claims that the risk of the adaptive estimate is not larger in order than $v_{\mathbb{I}}$ for all $\mathbb{I}$ such that $\delta_{\mathbb{I}}$ is sufficiently small.

THEOREM 4.1. *Let (3.5) hold true. Let an interval $\mathbb{I}$ be such that, for some $D \geq 0$, it holds with positive probability $\delta_{\mathbb{I}} \leq D$. Then*

$$\begin{aligned}&\mathbf{P}(\mathbb{I}\ \text{is rejected},\ \delta_{\mathbb{I}} \leq D)\\ (4.1)\quad &\qquad\leq \sum_{I\in\mathcal{I}(\mathbb{I})}\sum_{J\in\mathcal{J}(I)} 12\sqrt{e}\lambda_J(1+\log B)e^{-(\lambda_J-D)^2/(2a_\gamma)},\end{aligned}$$

*where $\lambda_J = \lambda(1 - s_\gamma \lambda N_J^{-1/2})$ with $N_J = \min\{|J|, |I\setminus J|\}$.*

*Moreover, if $N_J \geq 2s_\gamma \lambda$ for all $J \in \mathcal{J}(I)$ and all $I \in \mathcal{I}$, then it holds for the adaptive estimate $\hat{\theta} = \tilde{\theta}_{\hat{I}}$ on the random set $A = \{\mathbb{I}\ \text{is not rejected},\ \delta_{\mathbb{I}} \leq D\}$:*

$$|\tilde{\theta}_I - \tilde{\theta}_{\mathbb{I}}| \leq 2\lambda \tilde{v}_{\mathbb{I}}$$

*and*

$$|\tilde{\theta}_I - \theta_\tau| \leq (D + 3\lambda + 2\lambda s_\gamma(D+\lambda)|\mathbb{I}|^{-1/2})v_{\mathbb{I}}.$$

REMARK 4.1. It is easy to see that the sum on the right-hand side of the bound (4.1) can be made arbitrarily small by proper choice of the constant $\lambda$ and the sets $\mathcal{J}(I)$. Hence, the result of the theorem claims that with a dominating probability a "good" interval $\mathbb{I}$ will not be rejected and the adaptive estimate $\hat{\theta}$ is up to a constant factor as good as any of the "good" estimates $\tilde{\theta}_{\mathbb{I}}$.

REMARK 4.2. As mentioned in Remark 3.1, the probability bound on the right-hand side of (4.1) can be refined for the special case when the process $\theta_t$ is constant within $\mathbb{I}$ by replacing the factor $12\sqrt{e}\lambda_J(1+\log B)e^{-(\lambda_J-D)^2/(2a_\gamma)}$ by $6e^{-\lambda_J^2/(2a_\gamma)}$.

**5. Change point model.** A *change point* model is described by a sequence $T_1 < T_2 < \cdots$ of stopping times with respect to the filtration $\mathcal{F}_t$ and by values $\sigma_1, \sigma_2, \ldots$, where each $\sigma_k$ is $\mathcal{F}_{T_k}$-measurable. By definition, $\sigma_t = \sigma_k$ for $T_k < t \leq T_{k+1}$ and $\sigma_t$ is constant for $t < T_1$. This is an important special case of the model (2.1). For this special case the above procedure has a very natural interpretation: when estimating at the point $\tau$ we search for a largest interval of the form $[\tau - m, \tau]$ that does not contain a change point. This is done via testing for a change point within the candidate interval $I = [\tau - m, \tau]$. Note that the classical maximum likelihood test for no change point in the regression case with Gaussian $\mathcal{N}(0, \sigma^2)$ errors is also based on comparison of the mean values of observations $Y_t$ over the subintervals $I = [\tau - m, \tau - m']$ and every subinterval $J = [\tau - m', \tau]$ for different $m'$, so that the proposed procedure has strong appeal in this situation. However, there is an essential difference between testing for a change point and testing for homogeneity appearing as a building block of our adaptive procedure. Usually, a test for a change point is constructed in a way to provide the prescribed probability of a "false alarm," that is, rejecting the "no change point" hypothesis under homogeneity. Our adaptive procedure involves a lot of such tests for every candidate $I$, which leads to a multiple-testing problem. As a consequence, each particular test should be performed at a very high level; that is, it should be rather conservative providing a joint error probability at a reasonable level.

5.1. *Probability of a "false alarm."* For the change point model, a "false alarm" would mean that the candidate interval $I$ is rejected although the hypothesis of homogeneity is still fulfilled. The arguments used in the proof of Theorem 4.1 lead to the following upper bound for the probability of a "false alarm":

THEOREM 5.1. *If $I = [\tau - m, \tau]$ is an interval of homogeneity, that is, $\theta_t = \theta_\tau$ for all $t \in I$, then*

$$\mathbf{P}(I \text{ is rejected}) \leq \sum_{I \in \mathcal{I}(\mathbb{I})} \sum_{J \in \mathcal{J}(I)} 6 \exp\left(-\frac{\lambda^2}{2a_\gamma(1 + \lambda s_\gamma |J|^{-1/2})^2}\right).$$

This result is a special case of Theorem 4.1 with $\Delta_J \equiv 0$ when taking into account Remark 4.2.

Theorem 4.1 implies that for every fixed value $M$ there exists a fixed $\lambda$ providing a prescribed upper bound $\alpha$ for the probability of a "false alarm" for a homogeneous interval $I$ of length $M$. Namely, the choice

$$\lambda \geq (1 + \varepsilon)\sqrt{2a_\gamma \log \frac{M}{m_0 \alpha}}$$

leads for a proper small positive constant $\varepsilon > 0$ to the inequality

$$\sum_{I \in \mathcal{I}(\mathbb{I})} \sum_{J \in \mathcal{J}(I)} 6 \exp\left(-\frac{\lambda^2}{2a_\gamma(1 + \lambda s_\gamma |J|^{-1/2})^2}\right) \leq \alpha.$$

Here, $M/m_0$ is approximately the number of intervals in $\mathcal{J}(I)$ (see Section 6.1). This bound is, however, very rough, and it is only of theoretical importance since we estimate the probability of the sum of dependent events by the sum of single probabilities. The value of $\lambda$ providing a prescribed probability of a "false alarm" can be found by Monte Carlo simulation for the homogeneous model with constant volatility as described in Section 6.

5.2. *Sensitivity to change points and the mean delay.* The quality (sensitivity) of a change point procedure is usually measured by the mean delay between the occurrence of a change point and its detection.

To study this property of the proposed method, we consider the case of estimation at a point $\tau$ immediately after a change point $T_{\text{cp}}$. It is convenient to suppose that $T_{\text{cp}}$ belongs to the end points of an interval which is tested for homogeneity. In this case the "ideal" choice $\mathbb{I}$ is clearly $[T_{\text{cp}}, \tau]$. Theorem 4.1 claims that the quality of estimation at $\tau$ is essentially the same as if we knew the latest change point $T_{\text{cp}}$ a priori. In fact, one can state a slightly stronger assertion: every interval $I$ which is essentially larger than $\mathbb{I}$ will be rejected with high probability provided that the magnitude of the change is large enough.

Denote $m' = |\mathbb{I}|$, that is, $m' = \tau - T_{\text{cp}}$. Let also $I = [T_{\text{cp}} - m, \tau] = [\tau - m' - m, \tau]$ for some $m$, so that $|I| = m + m'$, and let $\theta$ (resp. $\theta'$) denote the value of the parameter $\theta_t$ before (resp. after) the change point $T_{\text{cp}}$. The magnitude of the change point is measured by the relative change $b = 2|\theta' - \theta|/\theta$.

It is worth mentioning that the values $\theta_t$ and especially $\theta_t'$ can be random and dependent on past observations. For instance, $\theta_t'$ may depend on $Y_t$ for all $t < T_{\mathrm{cp}}$.

The interval $I$ will certainly be rejected if $|\tilde{\theta}_{I \setminus \mathbb{I}} - \tilde{\theta}_{\mathbb{I}}|$ is sufficiently large compared to the corresponding critical value.

THEOREM 5.2. *Let $\mathbf{E}(Y_t | \mathcal{F}_{t-1}) = \theta$ before the change point at $T_{\mathrm{cp}}$ and $\mathbf{E}(Y_t | \mathcal{F}_{t-1}) = \theta'$ after it, and let $b = |\theta' - \theta|/\theta$. Let $I = [\tau - m' - m, \tau]$ with $m' = \tau - T_{\mathrm{cp}}$. If $\rho := \lambda s_\gamma / \sqrt{\min\{m, m'\}} < 1$ and*

$$(5.1) \qquad b \geq \frac{2\rho + \sqrt{2}\rho(1+\rho)}{1-\rho},$$

*then $\mathbf{P}(I \text{ is not rejected}) \leq 4e^{-\lambda^2/(2a_\gamma)}$.*

The result of Theorem 5.2 delivers some additional information about the sensitivity of the proposed procedure to change points. One possible question is about the minimal delay $m'$ between the change point $T_{\mathrm{cp}}$ and the first moment $\tau$ when the procedure starts to indicate this change point by selecting an interval of type $\mathbb{I} = [T_{\mathrm{cp}}, \tau]$. Due to Theorem 5.2, the change will be "detected" with high probability if the value $\rho = \lambda s_\gamma / \sqrt{m'}$ fulfills (5.1). With fixed $b > 0$, condition (5.1) leads to $\rho \leq bC_0$ for some fixed constant $C_0$. The latter condition can be rewritten in the form $m' \geq b^{-2} \lambda^2 s_\gamma^2 / C_0^2$. We see that this lower bound for the required delay $m'$ is proportional to $b^{-2}$, where $b$ is the change point magnitude. It is also proportional to the threshold $\lambda$ squared. In turn, for the prescribed probability $\alpha$ of rejecting a homogeneous interval of length $M$, the threshold $\lambda$ can be bounded by $C\sqrt{\log(M/m_0\alpha)}$. In particular, if we fix the length $M$ and $\alpha$, then $m' = O(b^{-2})$. If we keep fixed the values $b$ and $M$ but aim to provide a very small probability of a "false alarm" by letting $\alpha$ go to 0, then $m' = O(\log \alpha^{-1})$. All these issues are in agreement with the theory of change point detection; see, for example, Pollak (1985) and Brodsky and Darkhovsky (1993).

**6. LAVE in practice.** The aim of this section is to give some hints concerning the choice of the testing intervals and the smoothing parameter $\lambda$ and to illustrate the performance of the LAVE procedure on simulated and real data. We consider the simplest homogeneous model and we study the stability of the procedure in such a situation. Then a change point model is analyzed and the sensitivity with respect to the jump magnitude is measured. Finally, LAVE is applied to a set of exchange rate data.

6.1. *Choice of the sets $\mathcal{I}$ and $\mathcal{J}(I)$.* The presented algorithm involves the sets of interval candidates $\mathcal{I}$ and of testing intervals $\mathcal{J}(I)$. The simplest proposal is based on the use of a regular time grid $t_1, t_2, \ldots$, with grid step $m_0 \in \mathbb{N}$, that is, $t_k = m_0 k$, $k = 1, 2, \ldots$. For a given time point $\tau$, the set $\mathcal{I}$ of interval candidates is defined in the following way:

$$\mathcal{I} = \{I_k = [t_k, \tau] : t_k \leq \tau - m_0,\, k = 1, 2 \ldots\}.$$

Next, for every interval $I_k$, we define the set $\mathcal{J}(I_k)$ of testing subintervals $J_{k'} \subset I_k$ such that $J_{k'} = [t_{k'}, \tau]$ for all $t_{k'} > t_k$ belonging to the grid. The homogeneity within $I_k$ is then tested by comparing the pairs of estimates $\tilde{\theta}_J$ and $\tilde{\theta}_{I_k \setminus J}$ for all $J \in \mathcal{J}(I_k)$.

In this construction the sets $\mathcal{I}$, $\mathcal{J}(I)$ are completely determined by the grid step $m_0$. The value of $m_0$ should be selected possibly small, because it represents the minimal delay before the LAVE algorithm can detect a change point. Nevertheless, $m_0$ should be sufficiently large to provide stability of the estimates $\tilde{v}_J$ and $\tilde{v}_{I \setminus J}$. For the simulation and the analysis of real data we use $m_0 = 10$, which represents a good compromise. However, small changes in this value, that is, $5 \leq m_0 \leq 20$, do not appear to have great influence on the estimation results.

6.2. *Choice of $\lambda$ and $\gamma$.* The selection of $\gamma$ and, in particular, $\lambda$ is more critical. Theorem 5.1 suggests that in the context of a change point model, a reasonable approach for selecting $\lambda$ is by providing a prescribed level $\alpha$ for rejecting a homogeneous interval $I$ of a given length $M$. This would clearly imply at most the same level $\alpha$ for rejecting a homogeneous interval of a smaller length. However, the value of $\lambda$ which can be derived with the help of Theorem 5.1 is rather conservative. A more accurate choice can be made by Monte Carlo simulation. We examine the procedure described in Section 3.2 with the sets of intervals $\mathcal{I}$ and $\mathcal{J}(I)$ on the regular grid with the fixed step $m_0 = 10$. A constant (and therefore also time homogeneous) model assumes that the parameter $\theta_t$ does not vary in time, that is, $\theta_t \equiv \theta$. It can easily be seen that the value $\theta$ has no influence on the procedure under time homogeneity. One can therefore suppose that $\theta = 1$ and the original model (2.1) is transformed into the regression model $Y_t = 1 + s_\gamma \zeta_t$ with constant trend and homogeneous variance $s_\gamma$. This model is completely described, and, therefore, one can determine by simulation the value of $\lambda$ for which an interval of time homogeneity of length $M$ is not rejected with a frequency of 95%.

The values of $\lambda$ are computed for $M = 40$ and 80 and for the power transformations $\gamma = 0.5$, 1.0 and 2.0. The results are shown in Table 1. Note that the values of $\lambda$ calibrated for $M = 80$ are necessarily larger and therefore more conservative than the values of $\lambda$ calibrated for $M = 40$.

TABLE 1
*The value of $\lambda$, which, for a given power transformation $\gamma$, provides the rejection of an interval of time homogeneity of length $M$ with a frequency of 5%*

| Smoothing parameter | | | | | |
|---|---|---|---|---|---|
| **$\gamma = 0.5$** | | **$\gamma = 1.0$** | | **$\gamma = 2.0$** | |
| $M = 80$ | $M = 40$ | $M = 80$ | $M = 40$ | $M = 80$ | $M = 40$ |
| $\lambda = 2.74$ | $\lambda = 2.40$ | $\lambda = 2.58$ | $\lambda = 2.24$ | $\lambda = 2.18$ | $\lambda = 1.86$ |

6.3. *Simulation results for the change point model.* We now evaluate the performance of the LAVE algorithm on simulated data. Two change point time series of length $T = 240$ are considered. The simulated data display two jumps of the same magnitude in opposite directions: $\sigma_t = \sigma$ for $t \in [1, 80]$ and $t \in [161, 240]$ and $\sigma_t = \sigma'$ for $t \in [81, 160]$, where $\sigma = 1$ and $\sigma' = 3$ and 5, respectively. For each model 500 realizations are generated, and the estimation is performed at each time point $t \in [t_0, 240]$, where $t_0$ is set equal to 20.

We compute the estimation error for each combination of $\gamma$ and $\lambda$ with the following criterion:

$$(6.1) \qquad \sum_{t=20}^{240} \sum_{\omega=1}^{500} \left( \frac{\hat{\sigma}_t - \sigma_t}{\sigma_t} \right)^2 (\omega),$$

where the index $\omega$ indicates the realizations of the change point model. We note that in (6.1) the quadratic error is divided by the true volatility so that the criterion does not depend on the scale of $\sigma_t$. The results shown in Table 2 are favorable to the choice of the smaller value of $\gamma$, confirming that the loss of efficiency caused by $\gamma < 2$ is offset by the greater normality of the errors. Figures 2 and 3 show the results of the estimation for the power transformation $\gamma = 0.5$ and the value of $\lambda$ calibrated for an interval of time homogeneity of length $M = 40$ and $M = 80$, respectively. The plots on the top display the true process (straight line), the empirical median among all estimates (thick dotted line) and the empirical quartiles among all estimates (thin dotted lines). The plots on the bottom similarly display the length of the interval of time homogeneity, which is minimal (resp. maximal) just after (resp. just before) a change point, and the median and the quartiles among all estimates.

The results are satisfactory. The volatility is estimated precisely and the change points are quickly detected. As expected, the behavior of the method within homogeneous regions is very stable. The delay in detecting a change point becomes smaller as the jump size grows. Taking a smaller $\lambda$ also results in a smaller delay and improves the quality of estimation after the change points. The results for other power transformations look very similar and therefore are not displayed.

TABLE 2
*Estimation errors for all the combinations of parameters $\gamma$ and $\lambda$*

| | **Estimation error** | | | | | |
|---|---|---|---|---|---|---|
| | $\gamma = 0.5$ | | $\gamma = 1.0$ | | $\gamma = 2.0$ | |
| **Parameter** | $\lambda = 2.74$ | $\lambda = 2.40$ | $\lambda = 2.58$ | $\lambda = 2.24$ | $\lambda = 2.18$ | $\lambda = 1.86$ |
| Small jump | 19,241.9 | 17,175.3 | 19,121.2 | 16,522.5 | 24,887.2 | 17,490.9 |
| Large jump | 46,616.2 | 43,282.5 | 51,363.9 | 46,706.4 | 68,730.7 | 55,706.3 |

6.4. *Estimation of exchange rate volatility.* We apply the LAVE procedure to a set of nine exchange rates, which are available from the web site http:// federalreserve.gov of the U.S. Federal Reserve. The data sets represent daily exchange rates of the U.S. dollar (USD) against the following currencies: Australian dollar (AUD), British pound (BPD), Canadian dollar (CAD), Danish krone (DKR), Japanese yen (JPY), Norwegian krone (NKR), New Zealand dollar (NZD), Swiss franc (SFR) and Swedish krona (SKR). The period under consideration goes from January 1, 1990, to April 7, 2000. See Table 3.

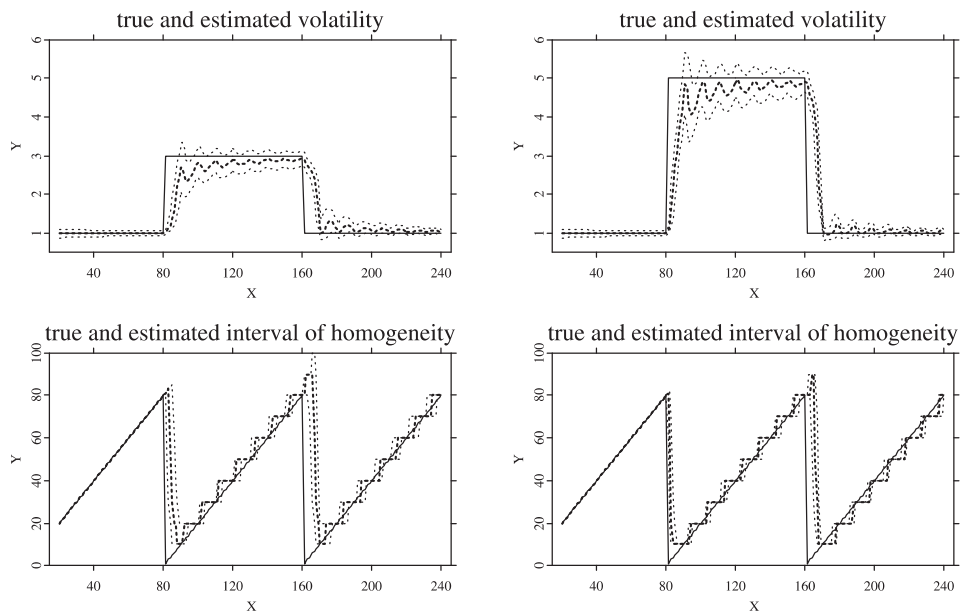

FIG. 2. *Estimation results for the change point model. The upper plots show the values of the standard deviation, while the lower plots show the values of the interval of homogeneity at each time point. True values* (solid line), *median of all estimates* (thick dotted line), *upper and lower quartiles* (thin dotted lines). *The value of $\lambda$ for $\gamma = 0.5$ and $M = 40$ has been used.*

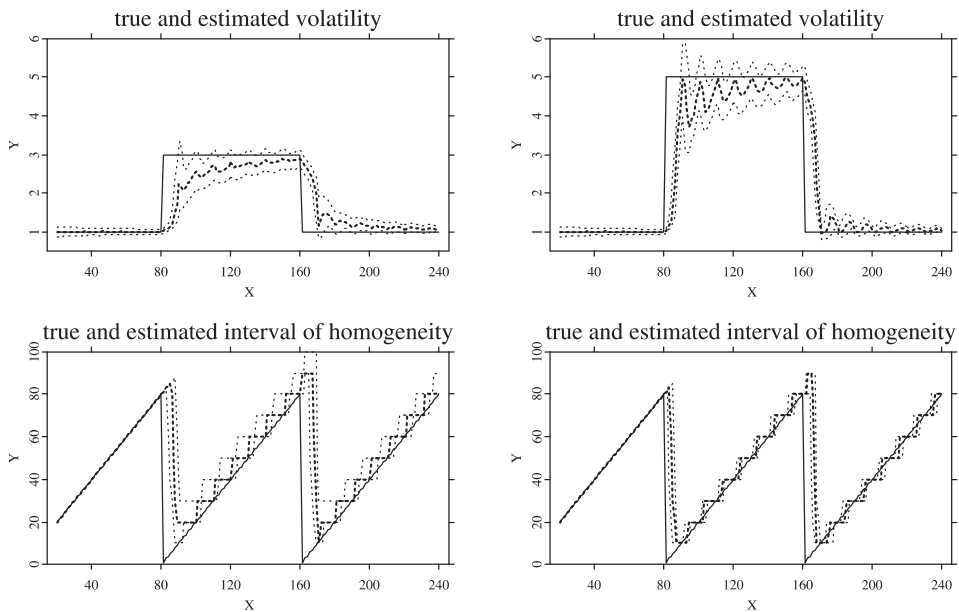

Fig. 3. *Estimation results for the jump model. The value of $\lambda$ for $\gamma = 0.5$ and $M = 80$ has been used.*

All the time series show qualitatively almost the same pattern; therefore, we provide the graphical example only for the two representative exchange rates JPY/USD and BPD/USD (Figure 4). The empirical mean of the returns is close to 0, while the empirical kurtosis is larger than 3. Fur-thermore, variance clustering and persistence of the autocorrelation of the square returns are also visible. The estimated standard deviation is nicely in accordance with the development of the volatility and, in particular, sharp changes in the volatility tend to be quickly recognized. Note also that the variability of the estimated interval of time homogeneity appears to grow as the estimated interval becomes larger. This is a feature of the algorithm because the number of tests grows with the accepted interval, so that a rejection becomes more probable. Nevertheless, this variability does not strongly affect the estimated volatility coefficient. Figure 5 shows the significantly persistent autocorrelation of the absolute returns, together with the autocorrelation of the absolute returns divided by the estimated standard deviation. The autocorrelation of the standardized absolute returns is not significant any more, and this fact supports the choice of a locally homogeneous model in order to explain the data.

Table 3
*Summary statistics*

| Currency | $n$ | Mean $\times 10^5$ | Variance $\times 10^5$ | Skewness | Kurtosis |
|---|---|---|---|---|---|
| AUD | 2583 | $-10.41$ | 3.191 | $-0.187$ | 8.854 |
| BPD | 2583 | $-0.679$ | 3.530 | $-0.279$ | 5.792 |
| CAD | 2583 | 8.819 | 0.895 | 0.042 | 5.499 |
| DKR | 2583 | 6.097 | 4.201 | $-0.037$ | 4.967 |
| JPY | 2583 | $-12.70$ | 5.486 | $-0.585$ | 7.366 |
| NKR | 2583 | 9.493 | 4.251 | 0.313 | 8.630 |
| NZD | 2583 | $-6.581$ | 3.604 | $-0.356$ | 49.17 |
| SFR | 2583 | 1.480 | 5.402 | $-0.186$ | 4.526 |
| SKR | 2583 | 12.66 | 4.615 | 0.372 | 9.660 |

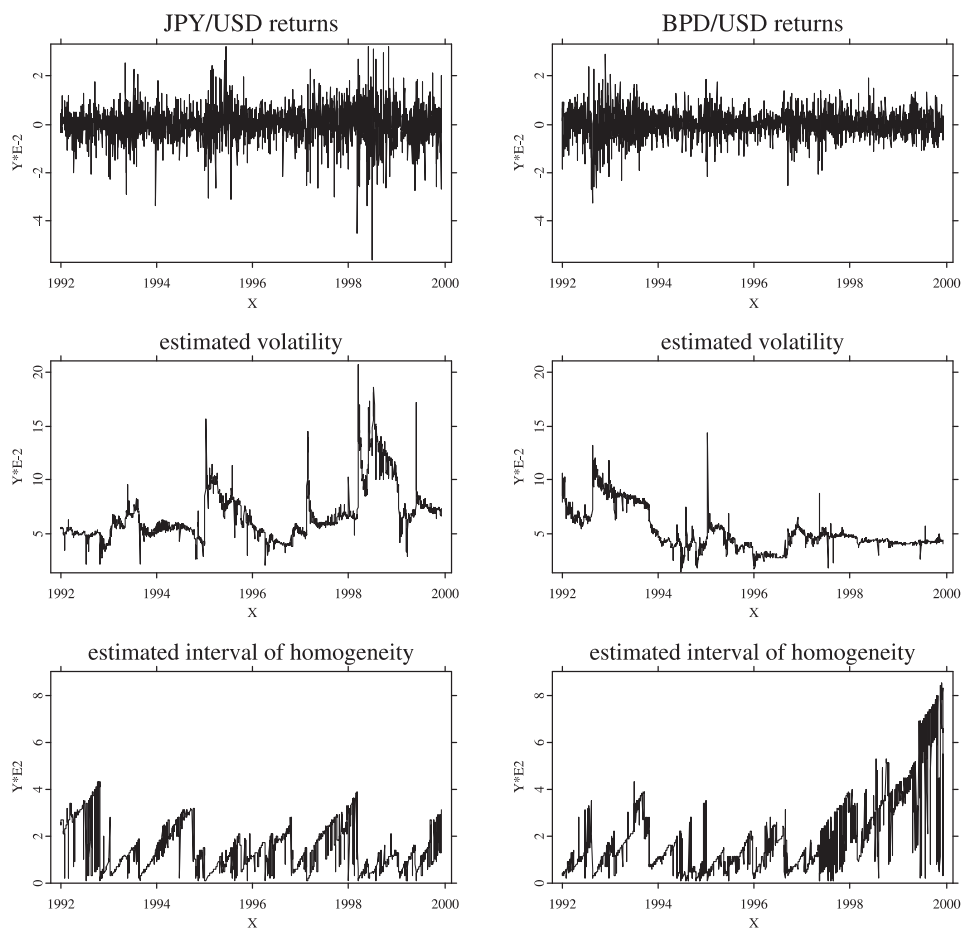

FIG. 4. *Exchange rate returns, estimated standard deviation and estimated interval of time homogeneity. The value of $\lambda$ for $\gamma = 0.5$ and $M = 80$ has been used.*

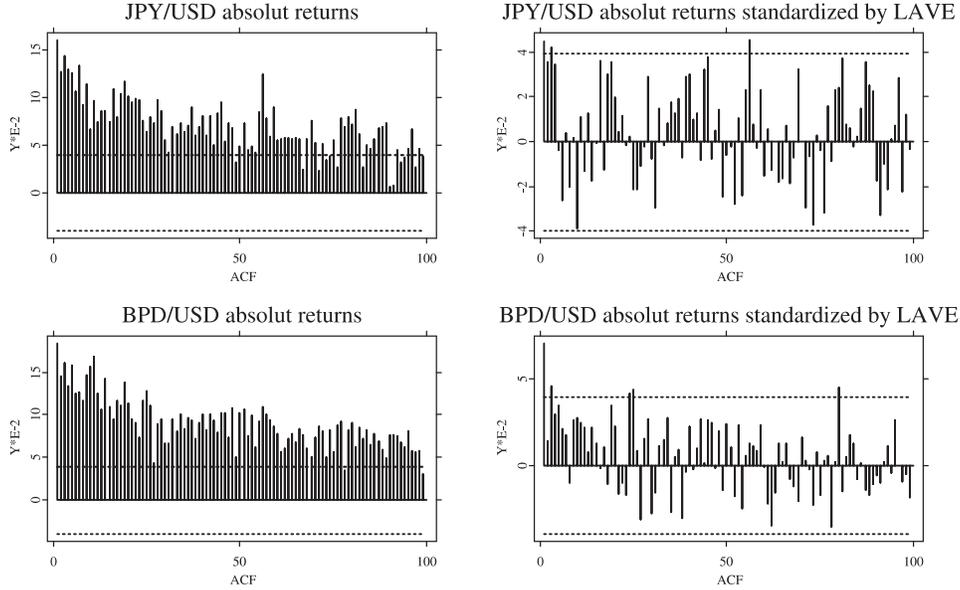

FIG. 5. *ACF of the absolute values of the exchange rate returns and ACF of the absolute values of the exchange rate returns standardized by LAVE.*

*A benchmark model.* As a matter of comparison, we also consider a model which is commonly used to estimate and forecast volatility processes: the GARCH$(1,1)$ model proposed by Bollerslev (1986):

$$\sigma_t^2 = \omega + \alpha R_{t-1}^2 + \beta \sigma_{t-1}^2.$$

Among all parametric volatility models, it represents the most common specification: "The GARCH$(1,1)$ is the leading generic model for almost all asset classes of returns. ... it is quite robust and does most of the work in almost all cases" [Engle (1995)].

We do not require the parameters to be constant throughout the whole sample, but, similarly to Franses and van Dijk (1996), we consider a rolling estimate. We thus fit the model to a sample of 350 observations, generate the forecast, delete the first observation from the sample and add the next one. Such a procedure reduces the harmful effect of possible parameter shifts on the forecasting performances of the model, even if at the same time it may increase the estimation variability.

The volatility is a hidden process which can be observed only together with a multiplicative error; therefore, the evaluation of the forecasting performance of an algorithm is not straightforward. Due to the model (2.1), $\mathbf{E}(R_{t+1}^2 | \mathcal{F}_t) = \sigma_{t+1}^2$. Therefore, given a forecast $\hat{\sigma}_{t+1|t}$, the empirical mean value of $|R_{t+1}^2 - \hat{\sigma}_{t+1|t}^2|^p$ can be used to measure the quality of this forecast. The forecast ability of the LAVE and the GARCH estimates is therefore evaluated with the following criterion:

$$\frac{1}{T - t_0 - 1} \sum_{t=t_0}^{T} |R_{t+1}^2 - \hat{\sigma}_{t+1|t}^2|^p \qquad \text{with } p = 0.5.$$

TABLE 4
*Forecast performance of LAVE relative to GARCH*

| Currency | $\gamma = 0.5$ | | $\gamma = 1.0$ | | $\gamma = 2.0$ | |
|---|---|---|---|---|---|---|
| | $M = 80$ | $M = 40$ | $M = 80$ | $M = 40$ | $M = 80$ | $M = 40$ |
| AUD | 0.942 | 0.945 | 0.963 | 0.962 | 0.991 | 0.982 |
| BPD | 0.961 | 0.960 | 0.979 | 0.970 | 1.006 | 1.000 |
| CAD | 0.974 | 0.979 | 0.989 | 0.992 | 1.010 | 0.997 |
| DKR | 0.978 | 0.980 | 0.985 | 0.987 | 1.010 | 1.004 |
| JPY | 0.951 | 0.949 | 0.971 | 0.966 | 1.006 | 0.997 |
| NKR | 0.961 | 0.957 | 0.972 | 0.965 | 0.998 | 0.984 |
| NZD | 0.878 | 0.879 | 0.904 | 0.902 | 0.952 | 0.947 |
| SFR | 0.985 | 0.984 | 0.992 | 0.990 | 1.004 | 1.000 |
| SKR | 0.965 | 0.961 | 0.973 | 0.968 | 0.982 | 0.977 |

The value of $p = 0.5$ is chosen instead of the more common $p = 2$ because we are interested in a robust criterion which is not too sensitive to the presence of outliers. The relative performance of the LAVE and the GARCH estimates is displayed in Table 4. The performance of the LAVE approach is clearly better; furthermore, the table gives a clear hint for the choice of the power transformation. Indeed, $\gamma = 0.5$ provides the smallest forecasting errors, while $\gamma = 2.0$ leads to the largest forecasting errors, which are sometimes larger than that of the GARCH model.

**7. Conclusions and outlook.** The locally adaptive volatility estimate (LAVE) is described and analyzed in this paper. It provides a nonparametric way for estimating and short-term forecasting the volatility of financial returns.

It is assumed that a local constant approximation of the volatility process holds over some unknown interval. The issue of filtering this interval of time homogeneity out of the return time series is considered, and a nonparametric approach is presented. The estimate of the volatility process is then found by averaging over the interval of time homogeneity.

A theoretical analysis of the properties of the LAVE algorithm is provided and the problem of selecting the smoothing parameters is analyzed through Monte Carlo simulation. The estimation results on change point models show that the method has reasonable performance in practice. An empirical application to exchange rate returns and a comparison with a GARCH(1, 1) also provide good evidence that the new method is competitive and can even outperform the standard parametric models, especially for forecasting with a short horizon.

An important feature of the proposed method is that it allows for a straightforward extension to multivariate volatility estimation; see Härdle, Herwartz and Spokoiny (2000) for a detailed discussion.

Obviously, if the underlying conditional distribution is not normal, the estimated volatility can give only partial information about the riskiness of the asset. Recent developments in risk analysis tend to focus on the estimation of the quantiles of the distribution. In this direction, the LAVE procedure can be used as a convenient tool for prewhitening the returns and obtaining a sample of "almost" identical and independently distributed returns, which do not display any more variance clustering. Therefore, the usual techniques of quantile estimation could be applied in a static framework. We regard such a development as a topic for future research.

**8. Proofs.** In this section, we collect the proofs of the results stated above. We begin by considering some useful properties of the power transformation introduced in Section 2.1.

*Some properties of the power transformation.* Let $g_\gamma(u)$ be the moment generating function of $\zeta_\gamma = D_\gamma^{-1}(|\xi|^\gamma - C_\gamma)$:

$$g_\gamma(u) = \mathbf{E} e^{u\zeta_\gamma}.$$

It is easy to see that this function is finite for $\gamma < 2$ and all $u$ and for $\gamma = 2$ and $u < 1$. For $\gamma = 1/2$, the function $2u^{-2} \log g_\gamma(u)$ is plotted in Figure 6.

LEMMA 8.1. *For every $\gamma \leq 1$ there exists a constant $a_\gamma > 0$ such that*

(8.1) $$\log \mathbf{E} e^{u\zeta_\gamma} \leq \frac{a_\gamma u^2}{2}.$$

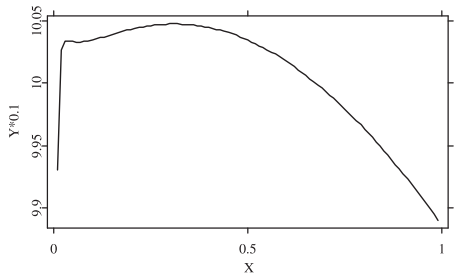

FIG. 6. *The log-Laplace transform of $\zeta_{1/2}$ divided by the log-Laplace transform of a standard normal r.v.*

PROOF. It is easy to check that the function $g_\gamma(u)$ with $\gamma \le 1$ is positive and smooth (infinitely many times differentiable). Moreover, the function $h_\gamma(u) = \log g_\gamma(u)$ is also smooth and satisfies $h_\gamma(0) = h'_\gamma(0) = 0$, $h''_\gamma(0) = \mathbf{E}\zeta_\gamma^2 = 1$. This yields that $u^{-2} h_\gamma(u) = u^{-2} \log g_\gamma(u)$ is bounded on every finite interval of the positive semiaxis $[0, \infty)$. It therefore remains to show that

$$\lim_{u \to \infty} u^{-2} \log \mathbf{E} e^{u\zeta_\gamma} < \infty.$$

Since $\zeta_\gamma(u) = D_\gamma^{-1}(|\xi|^\gamma - C_\gamma)$, it suffices to bound $u^{-2} \mathbf{E} e^{u|\xi|^\gamma / D_\gamma}$. For every $t > 0$,

$$\begin{aligned}
\mathbf{E} e^{u|\xi|^\gamma D_\gamma^{-1}} &= \mathbf{E} e^{u|\xi|^\gamma D_\gamma^{-1}} \mathbf{1}(|\xi| \le t) + \mathbf{E} e^{u|\xi|^\gamma D_\gamma^{-1}} \mathbf{1}(|\xi| > t) \\
&\le e^{ut^\gamma D_\gamma^{-1}} + \mathbf{E} e^{u|\xi| t^{\gamma-1} D_\gamma^{-1}} \\
&\le e^{ut^\gamma D_\gamma^{-1}} + 2 \mathbf{E} e^{u\xi t^{\gamma-1} D_\gamma^{-1}} \\
&= e^{ut^\gamma D_\gamma^{-1}} + 2 e^{u^2 t^{2\gamma-2} D_\gamma^{-2}}.
\end{aligned}$$

Next, with $t = u^{1/(2\gamma)}$ and $\gamma < 1$, for $u \to \infty$,

$$u^{-2} \log e^{ut^\gamma D_\gamma^{-1}} = u^{-1/2} D_\gamma^{-1} \to 0,$$

$$u^{-2} \log e^{u^2 t^{2\gamma-2} D_\gamma^{-2}} = u^{-(1-\gamma)/\gamma} D_\gamma^{-2} \to 0.$$

For $\gamma = 1$, the last expression remains bounded and the assertion follows. □

For $\gamma = 1/2$, condition (8.1) is satisfied with $a_\gamma = 1.005$.
The next technical statement is a direct consequence of Lemma 8.1.

LEMMA 8.2. *Let $c_t$ be a predictable process w.r.t. the filtration $\mathcal{F} = (\mathcal{F}_t)$; that is, every $c_t$ is a function of previous observations $R_1, \ldots, R_{t-1}$: $c_t = c_t(R_1, \ldots, R_{t-1})$. Then the process*

$$\mathcal{E}_t = \exp\left(\sum_{s=1}^t c_s \zeta_s - \frac{a_\gamma}{2} \sum_{s=1}^t c_s^2\right)$$

*is a supermartingale, that is,*

(8.2) $$\mathbf{E}(\mathcal{E}_t | \mathcal{F}_{t-1}) \le \mathcal{E}_{t-1}.$$

The next result has been stated in Liptser and Spokoiny (2000) for Gaussian martingales; however, the proof is based only on the property (8.2) and allows for straightforward extension to sums of the form $M_t = \sum_{s=1}^t c_s \zeta_s$.

THEOREM 8.1. *Let $M_t = \sum_{s=1}^{t} c_s \zeta_s$ with predictable coefficients $c_s$. Then let $T$ be fixed or a stopping time. For every $b > 0$, $B \geq 1$ and $\lambda \geq 1$,*

$$\mathbf{P}(|M_T| > \lambda \sqrt{\langle M \rangle_T}, b \leq \sqrt{\langle M \rangle_T} \leq bB) \leq 4\sqrt{e}\lambda(1 + \log B)e^{-\lambda^2/(2a_\gamma)},$$

*where*

$$\langle M \rangle_T = \sum_{t=1}^{T} c_t^2.$$

REMARK 8.1. If the coefficients $c_t$ are deterministic or independent of $M$, then Lemma 8.1 and the Chebyshev inequality yield

$$\mathbf{P}(|M_T| > \lambda \sqrt{\langle M \rangle_T}) \leq 2e^{-\lambda^2/(2a_\gamma)}.$$

PROOF OF THEOREM 3.1. Define

$$\bar{\theta}_I = \frac{1}{|I|} \sum_{t \in I} \theta_t, \qquad \xi_I = s_\gamma |I|^{-1} \sum_{t \in I} \theta_t \zeta_t.$$

Then $\tilde{\theta}_I = \bar{\theta}_I + \xi_I$. By the definition of $\Delta_I$,

(8.3) $$|\bar{\theta}_I - \theta_\tau| = |I|^{-1} \left| \sum_{t \in I} (\theta_t - \theta_\tau) \right| \leq \Delta_I.$$

Next, by (3.2)

$$\tilde{\theta}_I - \theta_\tau = \bar{\theta}_I - \theta_\tau + \xi_I,$$

and the use of (8.3) yields

$$\mathbf{P}(|\tilde{\theta}_I - \theta_\tau| > \Delta_I + \lambda v_I) \leq \mathbf{P}\left( \left| \sum_{t \in I} \theta_t \zeta_t \right| > \lambda \left( \sum_{t \in I} \theta_t^2 \right)^{1/2} \right).$$

In addition, if the volatility coefficient $\sigma_t$ satisfies $b \leq \sigma_t^2 \leq bB$ with some positive constants $b, B$, then the conditional variance $v_I^2 = s_\gamma^2 |I|^{-2} \sum_{t \in I} \theta_t^2$ satisfies

$$b'|I|^{-1} \leq v_I^2 \leq b'|I|^{-1}B,$$

with $b' = bs_\gamma^2$. Now the assertion follows from (3.5) and Theorem 8.1. □

PROOF OF THEOREM 3.2. It suffices to show that the inequalities $\Delta_I/v_I \leq D$ and

(8.4) $$|\tilde{\xi}_I| = |\tilde{\theta}_I - \bar{\theta}_I| \leq \lambda v_I$$

imply $|\tilde{\theta}_I - \theta_\tau| \leq \lambda' \tilde{v}_I$, where $\lambda'$ solves the equation $D + \lambda = \lambda'/(1 + \lambda' s_\gamma |I|^{-1/2})$. This would yield the desired result by Theorem 8.1; compare the proof of Theorem 3.1.

LEMMA 8.3. *Let $(\Delta_I/v_I)s_\gamma|I|^{-1/2} < 1$. Under (8.4)*

$$\tilde{v}_I \geq v_I\Big(\sqrt{1 - (\Delta_I/v_I)^2 s_\gamma^2 |I|^{-1}} - s_\gamma \lambda |I|^{-1/2}\Big) \geq v_I(1 - s_\gamma |I|^{-1/2}(\Delta_I/v_I + \lambda)).$$

PROOF. By the definition of $\tilde{v}_I$ in view of (8.4),

$$\tilde{v}_I = s_\gamma \tilde{\theta}_I |I|^{-1/2} \geq s_\gamma(\bar{\theta}_I - \lambda v_I)|I|^{-1/2}.$$

Since $\bar{\theta}_I$ is the arithmetic mean of $\theta_t$ over $I$,

$$\sum_{t \in I}(\theta_t - \bar{\theta}_I)^2 \leq \sum_{t \in I}(\theta_t - \theta_\tau)^2 \leq \Delta_I^2 |I|.$$

Next

$$s_\gamma^{-2}|I|v_I^2 = |I|^{-1}\sum_{t \in I}\theta_t^2 = \bar{\theta}_I^2 + |I|^{-1}\sum_{t \in I}(\theta_t - \bar{\theta}_I)^2 \leq \bar{\theta}_I^2 + \Delta_I^2,$$

so that

$$\bar{\theta}_I \geq s_\gamma^{-1}|I|^{1/2}v_I\sqrt{1 - (\Delta_I s_\gamma v_I^{-1}|I|^{-1/2})^2}.$$

Hence, under (8.4),

$$\tilde{v}_I \geq v_I\Big(\sqrt{1 - (\Delta_I s_\gamma v_I^{-1}|I|^{-1/2})^2} - s_\gamma \lambda |I|^{-1/2}\Big),$$

and the assertion follows. □

The bound (8.4) and the definition of $\Delta_I$ imply

$$|\tilde{\theta}_I - \theta_\tau| \leq |\bar{\theta}_I - \theta_\tau| + |\tilde{\theta}_I - \bar{\theta}_I| \leq \Delta_I + \lambda v_I \leq (D + \lambda)v_I.$$

By Lemma 8.3, $\tilde{v}_I \geq v_I(1 - s_\gamma D|I|^{-1/2} - s_\gamma \lambda |I|^{-1/2})$. Thus,

$$|\tilde{\theta}_I - \theta_\tau| \leq \frac{D + \lambda}{1 - s_\gamma(D + \lambda)|I|^{-1/2}}\tilde{v}_I = \lambda'\tilde{v}_I$$

as required. □

PROOF OF THEOREM 4.1. Let $\mathbb{I}$ be a "good" interval in the sense that, with high probability, $\Delta_J/v_J \leq D$ for some nonnegative constant $D$ and every $J \in \mathcal{J}(\mathbb{I})$. First we show that $\mathbb{I}$ will not be rejected with high probability provided that $\lambda$ is sufficiently large.

We proceed similarly as in the proofs of Theorems 3.1 and 3.2. The procedure involves the estimates $\tilde{\theta}_J$, $\tilde{\theta}_{I\setminus J}$ and the differences $\tilde{\theta}_J - \tilde{\theta}_{I\setminus J}$ for all $I \in \mathcal{I}(\mathbb{I})$ and all $J \in \mathcal{J}(I)$. The expansion $\tilde{\theta}_J = \bar{\theta}_J + \xi_J$ implies

$$\tilde{\theta}_J - \tilde{\theta}_{I\setminus J} = (\bar{\theta}_J - \bar{\theta}_{I\setminus J}) + (\xi_J - \xi_{I\setminus J}).$$

2224     D. MERCURIO AND V. SPOKOINYUnder the condition $\delta_I \leq D$,

$$|\bar{\theta}_J - \bar{\theta}_{I \setminus J}| \leq \Delta_I \leq D v_I \leq D \sqrt{v_J^2 + v_{I \setminus J}^2}.$$

Define the events

$$A_I = \bigcup_{J \in \mathcal{J}(I)} \left\{ |\xi_J - \xi_{I \setminus J}| \leq (\lambda_J - D)\sqrt{v_J^2 + v_{I \setminus J}^2} \right.$$

$$\left. \text{and } \sqrt{\frac{\tilde{v}_J^2 + \tilde{v}_{I \setminus J}^2}{v_J^2 + v_{I \setminus J}^2}} \geq 1 - s_\gamma \lambda N_J^{-1/2} \right\},$$

$$A_\mathbb{I} = \bigcup_{I \in \mathcal{I} : I \subseteq \mathbb{I}} A_I,$$

where $N_J = \min\{|J|, |I \setminus J|\}$ and $\lambda_J = \lambda(1 - s_\gamma \lambda N_J^{-1/2})$.

Define $A_\mathbb{I}^* = A_\mathbb{I} \cap \{\delta_\mathbb{I} \leq D\}$. On this set

$$\frac{|\tilde{\theta}_J - \tilde{\theta}_{I \setminus J}|}{\sqrt{\tilde{v}_J^2 + \tilde{v}_{I \setminus J}^2}} \leq \frac{|\bar{\theta}_J - \bar{\theta}_{I \setminus J}| + |\xi_J - \xi_{I \setminus J}|}{\sqrt{\tilde{v}_J^2 + \tilde{v}_{I \setminus J}^2}}$$

$$\leq (D + \lambda_J - D) \sqrt{\frac{v_J^2 + v_{I \setminus J}^2}{\tilde{v}_J^2 + \tilde{v}_{I \setminus J}^2}} \leq \frac{D + \lambda_J - D}{1 - s_\gamma \lambda N_J^{-1/2}} = \lambda.$$

It is easy to see that the conditional variance of $\xi_J - \xi_{I \setminus J}$ is equal to $v_J^2 + v_{I \setminus J}^2$. Arguing similarly to Lemma 8.3 and Theorem 3.1, we bound, with $\lambda_{J,D} = \lambda_J - D$,

$$\mathbf{P}(A_I) \leq \sum_{J \in \mathcal{J}(I)} \mathbf{P}\left( \frac{|\xi_J|}{v_J} > \lambda_{J,D} \right)$$

$$+ \mathbf{P}\left( \frac{|\xi_{I \setminus J}|}{v_{I \setminus J}} > \lambda_{J,D} \right) + \mathbf{P}\left( \frac{|\xi_J - \xi_{I \setminus J}|}{\sqrt{v_J^2 + v_{I \setminus J}^2}} > \lambda_{J,D} \right)$$

$$\leq \sum_{J \in \mathcal{J}(I)} 12 \sqrt{e} \lambda_J (1 + \log B) e^{-\lambda_{J,D}^2 / (2 a_\gamma)},$$

and the first assertion of the theorem follows.

Now we show that on the set $A_\mathbb{I}^*$ the estimate $\hat{\theta} = \tilde{\theta}_{\hat{I}}$ satisfies $|\hat{\theta} - \hat{\theta}_\mathbb{I}| \leq 2\lambda \tilde{v}_\mathbb{I}$.

Due to the above, on $A_\mathbb{I}^*$ the interval $\mathbb{I}$ will not be rejected and, hence $|\hat{I}| \geq |\mathbb{I}|$. Let $I$ be an arbitrary interval from $\mathcal{I}$ which is not rejected by the procedure. By construction $\mathbb{I}$ is one of the testing intervals for $I$. Denote

$J = I \setminus \mathbb{I}$. Note that $|I|(\tilde{\theta}_I - \tilde{\theta}_{\mathbb{I}}) = |J|(\tilde{\theta}_J - \tilde{\theta}_{\mathbb{I}})$, so that the event "$I$ is not rejected" implies $|\tilde{\theta}_J - \tilde{\theta}_{\mathbb{I}}| \leq \lambda\sqrt{\tilde{v}_J^2 + \tilde{v}_{\mathbb{I}}^2}$ and

$$|\tilde{\theta}_I - \tilde{\theta}_{\mathbb{I}}| \leq \frac{\lambda|J|}{|I|}\sqrt{\tilde{v}_J^2 + \tilde{v}_{\mathbb{I}}^2} \leq \frac{\lambda|J|}{|I|}(\tilde{v}_J + \tilde{v}_{\mathbb{I}}).$$

The use of $\tilde{v}_J = s_\gamma \tilde{\theta}_J |J|^{-1/2}$ and $|\tilde{\theta}_I - \tilde{\theta}_{\mathbb{I}}| \leq \lambda(\tilde{v}_J + \tilde{v}_{\mathbb{I}})$ yields

$$|\tilde{v}_J|J|^{1/2} - \tilde{v}_{\mathbb{I}}|\mathbb{I}|^{1/2}| \leq \lambda s_\gamma(\tilde{v}_J + \tilde{v}_{\mathbb{I}}),$$

implying

$$\tilde{v}_J \leq \frac{|\mathbb{I}|^{1/2} + \lambda s_\gamma}{|J|^{1/2} - \lambda s_\gamma}\tilde{v}_{\mathbb{I}}, \qquad \tilde{v}_J + \tilde{v}_{\mathbb{I}} \leq \frac{|J|^{1/2} + |\mathbb{I}|^{1/2}}{|J|^{1/2} - \lambda s_\gamma}\tilde{v}_{\mathbb{I}}.$$

Therefore,

$$|\tilde{\theta}_I - \tilde{\theta}_{\mathbb{I}}| \leq \frac{\lambda|J|(|J|^{1/2} + |\mathbb{I}|^{1/2})}{(|J| + |\mathbb{I}|)(|J|^{1/2} - \lambda s_\gamma)}\tilde{v}_{\mathbb{I}}.$$

It is straightforward to check that the function $f(x) = x^2(x+1)/[(x^2+1)(x-c)]$ with any $c \geq 0$ satisfies $f(x) \leq 2$ for all $x \geq 2c$. This implies with $x = |J|^{1/2}/|\mathbb{I}|^{1/2}$ and $c = \lambda s_\gamma/|\mathbb{I}|^{1/2}$ that

$$|\tilde{\theta}_I - \tilde{\theta}_{\mathbb{I}}| \leq 2\lambda \tilde{v}_{\mathbb{I}}$$

under the condition that $|J|^{1/2} \geq 2\lambda s_\gamma$.

Let $\Delta_{\mathbb{I}} \leq Dv_{\mathbb{I}}$. Similarly to Lemma 8.3, $\tilde{v}_{\mathbb{I}} \leq v_{\mathbb{I}}(1 + s_\gamma(D+\lambda)|\mathbb{I}|^{-1/2})$ and, by Theorem 3.1, $|\tilde{\theta}_{\mathbb{I}} - \theta_\tau| \leq (D+\lambda)v_{\mathbb{I}}$. This yields

$$|\tilde{\theta}_I - \tilde{\theta}_{\mathbb{I}}| \leq 2\lambda v_{\mathbb{I}}(1 + s_\gamma(D+\lambda)|\mathbb{I}|^{-1/2})$$

and

$$|\tilde{\theta}_I - \theta_\tau| \leq 2\lambda v_{\mathbb{I}}(1 + s_\gamma(D+\lambda)|\mathbb{I}|^{-1/2}) + (D+\lambda)v_{\mathbb{I}}$$
$$= (D + 3\lambda + 2\lambda s_\gamma(D+\lambda)|\mathbb{I}|^{-1/2})v_{\mathbb{I}}$$

as required. $\square$

PROOF OF THEOREM 5.2. To simplify the exposition we suppose that $\theta = 1$. (This does not restrict generality since one can always normalize each "observation" $Y_t$ by $\theta$.) We also suppose that $\theta' > 1$ and $b = \theta' - 1$. (The case when $\theta' < \theta$ can be considered similarly.) Finally, we assume that $m' = m$. (One can easily see that this case is the most difficult one.) We again apply the decomposition

$$\tilde{\theta}_J = 1 + \xi_J, \qquad \tilde{\theta}_{\mathbb{I}} = \theta' + \xi_{\mathbb{I}};$$

see the proof of Theorem 3.1. Hence

$$\tilde{\theta}_{\mathbb{I}} - \tilde{\theta}_J = b + \xi_{\mathbb{I}} - \xi_J.$$

It is straightforward to see that $v_J^2 = s_\gamma^2/m$ and $v_{\mathbb{I}}^2 = s_\gamma^2 \theta'/m$. By Lemma 8.1 (see also Remark 8.1)

$$\mathbf{P}(|\xi_J| > \lambda v_J) + \mathbf{P}(|\xi'| > \lambda v_{\mathbb{I}}) \le 4e^{-\lambda^2/(2a_\gamma)},$$

and it suffices to check that the inequalities $|\xi_J| \le \lambda v_J$, $|\xi_{\mathbb{I}}| \le \lambda v_{\mathbb{I}}$ and (5.1) imply

$$|\tilde{\theta}_J - \tilde{\theta}_{\mathbb{I}}| \ge \lambda\sqrt{\tilde{v}_J^2 + \tilde{v}_{\mathbb{I}}^2}.$$

Since $\theta' - 1 = b$ and since $\tilde{v}_J = s_\gamma |J|^{-1/2} \tilde{\theta}_J$ and similarly for $\tilde{v}_I$, we have under the conditions $|\xi_J| \le \lambda v_J$, $|\xi_{\mathbb{I}}| \le \lambda v_{\mathbb{I}}$,

$$|\tilde{\theta}_J - \tilde{\theta}_{\mathbb{I}}| \ge b - \frac{\lambda s_\gamma(\theta' + 1)}{\sqrt{m}} = b(1-\rho) - 2\rho,$$

$$\tilde{v}_J = \frac{s_\gamma}{\sqrt{m}}(1 + \xi_J) \le \lambda^{-1}\rho(1+\rho),$$

$$\tilde{v}_{\mathbb{I}} = \frac{s_\gamma}{\sqrt{m}}(1 + \xi_{\mathbb{I}}) \le \lambda^{-1}\rho(1+\rho),$$

with $\rho = m^{-1/2}\lambda s_\gamma$. Therefore,

$$|\tilde{\theta}_J - \tilde{\theta}_{\mathbb{I}}| - \lambda\sqrt{\tilde{v}_J^2 + \tilde{v}_{\mathbb{I}}^2} \ge b(1-\delta) - 2\rho - \sqrt{2}\rho(1+\rho) > 0$$

in view of (5.1), and the assertion follows. $\square$

## REFERENCES

BOLLERSLEV, T. (1986). Generalized autoregressive conditional heteroskedasticity. *J. Econometrics* **31** 307–327. MR853051

BRODSKY, B. and DARKHOVSKY, B. (1993). *Nonparametric Methods in Change-Point Problems.* Kluwer, Dordrecht.

CAI, Z., FAN, J. and LI, R. (2000). Efficient estimation and inferences for varying-coefficient models. *J. Amer. Statist. Assoc.* **95** 888–902. MR1804446

CAI, Z., FAN, J. and YAO, Q. (2000). Functional-coefficient regression models for nonlinear time series. *J. Amer. Statist. Assoc.* **95** 941–956. MR1804449

CARROLL, R. and RUPPERT, D. (1988). *Transformation and Weighting in Regression* Chap. 4. Chapman and Hall, New York. MR1014890

CLEMENTS, M. P. and HENDRY, D. F. (1998). *Forecasting Economic Time Series.* Cambridge Univ. Press. MR1662973

ENGLE, R. F., ed. (1995). *ARCH, Selected Readings.* Oxford Univ. Press.

ENGLE, R. F. and BOLLERSLEV, T. (1986). Modelling the persistence of conditional variances (with discussion). *Econometric Rev.* **5** 1–87. MR876792

Fan, J., Jiang, J., Zhang, C. and Zhou, Z. (2003). Time-dependent diffusion models for term structure dynamics. Statistical applications in financial econometrics. *Statist. Sinica* **13** 965–992. MR2026058

Fan, J. and Zhang, W. (1999). Statistical estimation in varying coefficient models. *Ann. Statist.* **27** 1491–1518. MR1742497

Franses, P. and van Dijk, D. (1996). Forecasting stock market volatility using (nonlinear) GARCH models. *J. Forecasting* **15** 229–235.

Glosten, L., Jagannathan, R. and Runkle, D. (1993). On the relation between the expected value and the volatility of the nominal excess return on stocks. *J. Finance* **48** 1779–1801. MR1242517

Gouriéroux, C. (1997). *ARCH Models and Financial Applications*. Springer, Berlin. MR1439744

Härdle, W., Herwartz, H. and Spokoiny, V. (2000). Time inhomogeneous multiple volatility modelling. Unpublished manuscript.

Härdle, W. and Stahl, G. (1999). Backtesting beyond var. Technical Report 105, Sonderforschungsbereich 373, Humboldt Univ., Berlin.

Harvey, A., Ruiz, E. and Shephard, N. (1994). Multivariate stochastic variance models. *Rev. Econom. Stud.* **61** 247–264.

Kitagawa, G. (1987). Non-Gaussian state-space modelling of nonstationary time series (with discussion). *J. Amer. Statist. Assoc.* **82** 1032–1063. MR922169

Lepski, O. (1990). On a problem of adaptive estimation in Gaussian white noise. *Theory Probab. Appl.* **35** 454–466. MR1091202

Lepski, O. and Spokoiny, V. (1997). Optimal pointwise adaptive methods in nonparametric estimation. *Ann. Statist.* **25** 2512–2546. MR1604408

Liptser, R. and Spokoiny, V. (2000). Deviation probability bound for martingales with applications to statistical estimation. *Statist. Probab. Lett.* **46** 347–357. MR1743992

Mikosch, T. and Starica, C. (2000). Change of structure in financial time series, long range dependence and the garch model. Technical Report 58, Aarhus School of Business, Aarhus Univ.

Nelson, D. B. (1991). Conditional heteroscedasticity in asset returns: A new approach. *Econometrica* **59** 347–370. MR1097532

Pollak, M. (1985). Optimal detection of a change in distribution. *Ann. Statist.* **13** 206–227. MR773162

Sentana, E. (1995). Quadratic ARCH models. *Rev. Econom. Stud.* **62** 639–661.

Spokoiny, V. (1998). Estimation of a function with discontinuities via local polynomial fit with an adaptive window choice. *Ann. Statist.* **26** 1356–1378. MR1647669

Taylor, S. J. (1986). *Modelling Financial Time Series*. Wiley, Chichester.

Humboldt-Universität zu Berlin
Spandauer Strasse 1
10178 Berlin
Germany
e-mail: isedoc02@wiwi.hu-berlin.de

Weierstrass-Institute
Mohrenstrasse 39
10117 Berlin
Germany
e-mail: spokoiny@wias-berlin.de